\documentclass[11pt,letterpaper]{amsart}
\usepackage{amsmath,amsfonts,amssymb,amscd,amsthm,amsbsy,epsf}
\usepackage[dvips]{graphicx}
\usepackage{comment}
\usepackage[dvips]{color}
\usepackage{xcolor}
\usepackage{orcidlink}
\usepackage{enumitem}

\numberwithin{equation}{section}
\newtheorem{theorem}{Theorem}
\newtheorem{meta-thm}[theorem]{Meta-Theorem}
\newtheorem{lemma}[theorem]{Lemma}

\theoremstyle{remark}
\newtheorem{remark}[theorem]{Remark}

\newcommand{\noaverage}[1]{({#1})^0}

\def\D{{\mathcal D}}

\def\LL{ {\mathcal{L}}}

\def\complex{{\mathbb C}}
\def\integer{{\mathbb Z}}
\def\nat{{\mathbb N}}

\def\real{{\mathbb R}}

\def\torus{{\mathbb T}}

\def\eps{\varepsilon}

\def\th{\theta}

\def\<{\langle}
\def\>{\rangle}
\def\froeschle{Froeschl\'e}

\begin{document}
\title[Domains of analyticity of lower dimensional tori]{Computation of Domains of analyticity of lower dimensional tori in a weakly dissipative Froeschlé map }

\author[A.P. Bustamante]{Adri\'an P. Bustamante}  
\address{
Department of Mathematics and Physics, University of Roma Tre, Via della Vasca Navale 84, 00146 Roma (Italy)}
\email{adrian.perezbustamante@uniroma3.it}
\thanks{The author has been supported by the research project  MIUR-PRIN 2020XBFL “Hamiltonian and Dispersive PDEs”}



\begin{abstract}
  We consider a Froeschl\'e map and we add a weak dissipation of the form $\lambda(\eps) = 1- \eps^3$, where $\eps$ is the parameter of perturbation. We compute formal expansions of lower dimensional tori, both in the conservative and weakly dissipatives cases, and use them to estimate the shape of their domains of analyticity as functions of $\eps$. Our results support conjectures in the literature.

 \end{abstract} 
\maketitle


\section{Introduction}\label{sec:intro}

We consider a classical Froeshlé map \eqref{eq:map} with a trigonometric perturbation and a weak dissipation of the form $\lambda(\eps) =1-\eps^3$, where $\eps$ is the parameter of perturbation. The choice of this type type of dissipation is motivated by the fact that some problems of physical interest can be described by adding a small dissipation to a conservative system, e.g., Hamiltonian systems with dissipation proportional to the velocity -like in some problems appearing in celestial dynamics, see \cite{Cel-10}.

In this work we compute formal expansions of lower dimensional tori, for the map \eqref{eq:map}, and use them to approximate their domains of analyticity as functions of a complex parameter $\eps\in \complex$. We recall  that the theoretical results in \cite{Jo-Lla-Zou-99,GallavottiG02, Bus-Lla-23} suggest the power series are not convergent and that, in fact, are Gevrey. One of the main objectives to approximate numerically the analyticity domains is to explore some conjectures in the literature \cite{Jo-Lla-Zou-99,GallavottiG02} about the optimal domain of analyticity of the formal expasions in the conservative case. 

As a byproduct of our computation we also explore the Gevrey character of the Lindstedt expansions. We recall that in \cite{Bus-Lla-23} was proved that the expansions of lower dimensional tori for a family of Froeschlé maps are Gevrey. The estimates of the Gevrey exponent in \cite{Bus-Lla-23} are very simple and are not expected to be optimal, so in this work we also perform a numerical approximation of the Gevrey exponents. 

Our approach is to compute numerically the Lindsted series of lower dimensional tori, in the conservative and dissipative case, by following the approach in \cite{Bus-Lla-23}. Once the Linstedt series are computed up to very high order we compute the poles of Padé and Log-Padé approximants (see Appendix \ref{sec:appendix}) to approximate the domains of analyticity of lower dimensional tori as functions of $\eps\in \complex$. The use of Lindstedt series and Padé approximants has been used systematically in the literature to study breakdown thresholds and domains of analyticity \cite{Ber-Chi-90, Ber-Cel-Chi-Fal-92,Fal-Lla-92, BustamanteC19, Bus-Cel-Lho-23}. 

The paper is organized as follows. In Section \ref{sec:prelims} we describe the settings of the problem and give some definitions. Section \ref{sec:lindstedt} briefly describes the algorithm to compute the Lindstedt series both in the conservative and dissipative case. Finally, in Section \ref{sec:num-results} we present the results of our computations. We have included Appendix \ref{sec:appendix} describing the computation of Padé and Log-Padé approximants.

\section{Preliminaries}\label{sec:prelims}
We consider a \froeschle-like map defined on the cylinder $\torus^2\times \real^2$ given by \begin{align}
    p_{n+1} &= (1-\gamma\eps^3)p_n + \mu_\eps + \eps \nabla V(q_n) \nonumber\\
    q_{n+1} &= q_n + p_{n+1}\label{eq:map}
\end{align}
where $q_n \in\torus^2$, $p_n\in\real^2$, $V(x_1, x_2) := -\cos(x_1)-\cos(x_2) - \cos(x_1 + x_2)$, and $\eps$ is a parameter of perturbation. We use a parameter $\gamma\in\{0,1\}$ to distinguish between the classical \froeschle\  map ($\gamma = 0$) and a dissipative version of the map ($\gamma=1$). 

Note that when $\gamma=1$ we have a singular perturbation due to the fact that the perturbation also introduces a \textit{weakly} dissipative term $\lambda(\eps) = 1-\eps^3$. Since the addition of a dissipative term affects drastically the existence of quasi-periodic orbits an extra parameter, $\mu_\eps\in\real^2$, needs to be adjusted to ensure the existence of orbits of a given frequency $\omega\in\real$. When $\gamma=0$ the map remains symplectic and the extra parameter is not needed.

\subsection{Quasi-periodic orbits}
We are interested in computing formal expansions of lower dimensional tori for the map \eqref{eq:map}, so we consider a frequency $\omega\in\real$ satisfying the Diophantine condition \begin{equation}\label{eq:diophantine}
|m \omega -n|\geq \nu |m|^{-\ell} \qquad n\in\integer, \, m\in \integer\backslash\{0\}
\end{equation}
for some  $\nu, \ell > 0 $. We say that  $\omega\in \D(\ell)$ if $\omega\in\real$ satisfy \eqref{eq:diophantine}.   

For our purpose it is convenient to write the map \eqref{eq:map} as a second order equation for orbits in half of the variables. Using the second equation in \eqref{eq:map} one has $p_{n+1} =  q_{n+1} - q_n$, substituting this last expression into the first equation of $\eqref{eq:map}$ one obtains \begin{equation}\label{eq:second-order-map}
q_{n+1} - (2 -\gamma\eps^3)q_n + (1-\gamma\eps^3)q_{n-1}= \eps\nabla V(q_n) +\mu_\eps.
\end{equation} 

Following \cite{Bus-Lla-23, Jo-Lla-Zou-99} we recall that quasi-periodic orbits of a given frequency, for the map \eqref{eq:second-order-map}, can be found using a parametric representation of the variable $q_n\in \torus^2$. We say that a sequence $q_n$ is quasiperiodic, of frequency $\omega\in\real$, if there exists a function, $h$, such that \begin{equation}
    q_n=h(\omega n )
\end{equation}
where $h:\torus\longrightarrow\torus^2$ is  oftentimes called \textit{hull} function. One can think of $h$ as an embedding of $\torus$ into the phase space. 

We recall that the hull function must satisfy the condition  \begin{equation}\label{eq:top-hull}
    h(\theta +2\pi ) = h(\theta) + 2\pi k
\end{equation} 
where $k=(k_1, k_2)\in \integer^2$, which is a topological property of the embedding. The geometric meaning of \eqref{eq:top-hull} is that $h$ gives a parameterization of a circle embedded in $\torus^2$ which winds $k_1$ times in the first component and $k_2$ times in the second one. We also note that the hull function can be written in terms of a periodic function, $g:\real\longrightarrow\real^2$, such that \begin{equation}\label{eq:periodic-part}
    h(\theta) = \theta k + g(\theta). 
\end{equation}
We call $g$ the periodic part of the hull function. 

\subsection{The invariance equation and normalization.}
If $\omega$ is irrational, finding a quasiperiodic orbit $\{q_n\}$ of \eqref{eq:second-order-map} is equivalent to finding a hull function, and a drift parameter $\mu_\eps$, satisfying \begin{equation}\label{eq:inv-hull}
    h(\theta + \omega) - (2-\gamma \eps^3)h(\theta) + (1-\gamma \eps^3) h(\theta-\omega) =\eps \nabla V(h(\theta)) + \mu_\eps.
\end{equation} 
Note that equation \eqref{eq:inv-hull} is underdetermined. That is, if $h$ is a solution of \eqref{eq:inv-hull} then, for any $\sigma\in \torus$ the function \begin{equation}\label{eq:under-hull}h_\sigma(\theta) := h(\theta+\sigma)\end{equation} is also a solution of equation \eqref{eq:inv-hull}. The geometric meaning of this underdeterminancy is the choice of the origin of the system of coordinates of the reference manifold $\torus$. It is important to choose a normalization that fixes this underdeterminancy in order to obtain a unique solution of \eqref{eq:inv-hull}, we follow \cite{Bus-Lla-23, Jo-Lla-Zou-99} and choose a normalization condition for the periodic part of the hull function. 

Using \eqref{eq:periodic-part} in equation \eqref{eq:inv-hull}  one obtains the following invariance equation for the periodic part of  the hull function \begin{equation}
    \label{eq:inv-perpart} g(\theta +\omega) -2g(\theta) + g(\theta-\omega) + \gamma\eps^3\left(k\omega +g(\theta)-g(\theta-\omega)\right) = \eps \nabla V(\theta k + g(\theta)) +\mu_\eps.
\end{equation} 
In terms of the function $g$ the underdeterminancy \eqref{eq:under-hull} translates into the fact that if $g$ is a solution of \eqref{eq:inv-perpart} then, for any $\sigma\in \real$, the function \begin{equation}\label{eq:under-per} 
g_\sigma(\theta) := \sigma k + g(\theta + \sigma)
\end{equation}
is also a solution of \eqref{eq:inv-perpart}.
To deal with this underdeterminancy we consider the following normalization
\begin{equation}\label{eq:normalization}
    \int_\torus k\cdot g(\theta)d\theta =0. 
\end{equation}
In \cite{Bus-Lla-23} it was shown that with the normalization \eqref{eq:normalization} the solution of equation \eqref{eq:inv-perpart}
is unique.

\section{Lindstedt series}\label{sec:lindstedt}
In this section we briefly describe a procedure to compute formal solutions \begin{equation}\label{eq:g-mu-expan}
g_\eps(\theta) = \sum_{n=0}^\infty g_n(\theta)\eps^n\qquad \mu_\eps = \sum_{n=0}^\infty \mu_n \eps^n
\end{equation} of \eqref{eq:inv-perpart}. We note that the formal series $g_\eps$ in \eqref{eq:g-mu-expan} satisfy the normalization \eqref{eq:normalization} if, for any $n\geq 0$, \begin{equation}
    \int_\torus k\cdot g_n(\theta) d\theta = 0. \label{eq:ord-n-normalization}
\end{equation}
A very detailed proof of the existence of the Lindstedt series, and the uniqueness under normalization, can be found in \cite{Bus-Lla-23}. In the following paragraphs we summarized the steps needed to compute the formal expansions. 

\subsection{Formal expansion of the potential $\nabla V$}
To find the expansions \eqref{eq:g-mu-expan} the right hand side of \eqref{eq:inv-perpart} also needs to be expanded in power series, that is,  one needs to find the expansion
\begin{equation} \label{eq:V-expansion}
    \eps\nabla V(\theta k + g_\eps(\theta)) = \sum_{n=1}^\infty R_n(\theta)\eps^n.
\end{equation}
To compute the coefficients $R_n$ in \eqref{eq:V-expansion} we use an idea that has been used many times in automatic differentiation \cite{BrentK78}. Assume that $V:\real^2\longrightarrow\real$ is a trigonometric polynomial of degree $J$,  then \begin{equation}\label{eq:trig-pol} 
\nabla V(q) = \sum_{0<|\ell|\leq J} \alpha_\ell e^{i k\cdot q}.
\end{equation}
Thus, to find the expansion \eqref{eq:V-expansion} it is enough to compute expansions for the functions $e^{i\ell\cdot (\theta k + g_\eps(\theta) )}$. We write \begin{equation}\label{eq:exp-expansion}
    e^{i\ell\cdot (\theta k + g_\eps(\theta) )} = \sum_{n=0}^\infty F^{\ell,k}_n(\theta)\eps^n
\end{equation}
and note that \begin{equation}\label{eq:deriv-eps}
    \frac{d}{d\eps}e^{i\ell\cdot (\theta k + g_\eps(\theta) )} = i\ell \cdot \frac{d}{d\eps}g_\eps(\theta)e^{i\ell\cdot (\theta k + g_\eps(\theta) )}. 
\end{equation}
Using the formal expansions \eqref{eq:g-mu-expan} and \eqref{eq:exp-expansion} in \eqref{eq:deriv-eps}, we see that the coefficients $F_n^{\ell, k}$ in \eqref{eq:exp-expansion} must satisfy  the following relations \begin{equation}
    F_{n}^{\ell, k}(\theta)= \frac{i}{n} \sum_{m=0}^{n-1} (m+1)\ell\cdot g_{m+1}(\theta) F_{n-1-m}^{\ell,k}(\theta)
\end{equation}
for $n\geq 1$ and $F_0^{\ell,k}(\theta)= e^{i\ell\cdot (\theta k + g_0)}$. Finally, for $n\geq 1$, \begin{equation}
    R_n(\theta) =\sum_{0<|\ell|\leq J} \alpha_\ell F_{n-1}^{\ell,k}(\theta).\label{eq:R-n}
\end{equation}
Note that $R_n$ depends on $g_0, g_1, \dots, g_n$.

\subsection{Computation of the Lindstedt series} 

The computation of the formal power series \eqref{eq:g-mu-expan} solving \eqref{eq:inv-perpart} is different whether we are in the conservative case ($\gamma = 0$) or in the weakly dissipative case ($\gamma\neq 0$).

\subsubsection{Conservative case.}
In this case we have $\gamma = 0$ and the drift parameter $\mu_\eps \equiv 0$. Plugging \eqref{eq:g-mu-expan} and \eqref{eq:V-expansion} into \eqref{eq:inv-perpart} one has that at each order one has to solve the cohomology equations \begin{equation}\label{eq:g_n}
    g_n(\theta + \omega) -2g_n(\theta) +g_n(\theta -\omega) = R_n(\theta)
\end{equation}
where $R_n$ is given by \eqref{eq:R-n}. Functional equations of the same type as  \eqref{eq:g_n} are well known in KAM theory and its solubility is assured as long as the right hand size has zero average. For the sake of completeness and to be able to refer to it we recall the following result, see \cite{Llave01}. 
\begin{lemma}
    Let $\omega$ Diophantine and $\eta:\torus\longrightarrow \real^n$ be a continuous function with zero average, that is $\int_\torus \eta(\theta)d\theta = 0$. Then there exists an $L^2$ function, $\varphi$, solving the equation \begin{equation}
        \varphi(\theta + \omega) - 2\varphi(\theta) + \varphi(\theta-\omega) = \eta(\theta).
        \label{eq:cohom-eq}
    \end{equation}
    Moreover, the solution can be written as \begin{equation*}
        \varphi(\theta) = \hat\varphi_0 + \sum_{\ell \in \integer\backslash \{0\}}\frac{\hat{\eta}_\ell}{2 (\cos(\ell \omega) - 1)} e^{i\ell \theta}
    \end{equation*}
    where $\hat\eta_\ell$ are the Fourier coefficients of $\eta(\theta)$. Note that $\hat\varphi_0 $ is a free parameter. 
\end{lemma}

At order order $n=0$ one has the equation \begin{equation}\label{eq:ord-zero}
    g_0(\theta + \omega) - 2g_0(\theta) + g_0(\theta -\omega) = 0
\end{equation}
which is solved by any constant $g_0$. To find a solution satisfying the normalization \eqref{eq:ord-n-normalization} we choose a constant orthogonal to $k$, that is, we fix a vector $k^\perp$ orthogonal to $k$ and take $g_0 = \beta_0 k^\perp$, for some $\beta_0\in \real$ that will be fixed in the next step.

The equation at order $n=1$ is \begin{equation}\label{eq:eq-ord-1}
    g_1(\theta+\omega) - 2g_1(\theta) + g_1(\theta - \omega) = \nabla V(\theta k + \beta_0 k^\perp ),
\end{equation}
which has a solution if  \begin{equation}
    \int_\torus \nabla V(\theta k + \beta_0 k^\perp) d\theta = 0.\label{eq:zero-ave-ord-1} 
\end{equation}
We note that if $\nabla V$ is a real valued trigonometric polynomial then, there exist at least two values of $\beta_0$  for which \eqref{eq:zero-ave-ord-1} is satisfied, see \cite{Bus-Lla-23}. The selection of $\beta_0$ satisfying \eqref{eq:zero-ave-ord-1} completely determines $g_0$. 

Once that \eqref{eq:zero-ave-ord-1} is fulfilled we can find a solution of equation \eqref{eq:eq-ord-1} up to an additive constant, we chose this constant in the direction of $k^\perp$ to attain the normalization \eqref{eq:ord-n-normalization}. That is, we consider a solution of \eqref{eq:eq-ord-1} of the form $g_1 + \beta_1 k^\perp $, where $g_1$ satisfies \eqref{eq:eq-ord-1} and has zero average. The constant $\beta_1\in \real$ is chosen in the next step, so the average of  the right hand side of equation \eqref{eq:g_n} vanishes. This procedure will be repeated at each order $n$, we summarize it in the following Lemma.

\begin{lemma}[\cite{Jo-Lla-Zou-99, Bus-Lla-23}]\label{lem:solution-induction}
    Let $g_0 = \beta_0 k^\perp$ be the choice for the solution of equation \eqref{eq:ord-zero} satisfying \eqref{eq:zero-ave-ord-1}. Assume the non-degeneracy condition \begin{equation}
        \int_\torus k^\perp\cdot D^2V(\theta k + \beta_0 k^\perp) k^\perp d\theta \neq 0.\label{eq:non-degeneracy}
    \end{equation} is satisfied.
    For $n\geq 2$, if we have that $g_0, \dots, g_{n-1} + \beta_{n-1}k^\perp $ solve the invariance equation \eqref{eq:g_n} of order $0,\dots, n-1$, respectively, as well as the normalization \eqref{eq:ord-n-normalization} ($g_0, \dots, g_{n-1}$ are uniquely determined, $\beta_{n-1}$ is arbitrary), then:
    \begin{itemize}
        \item There is only one $\beta_{n-1}$ so that equation \eqref{eq:g_n}, at order $n$, has a solution.
        \item All the solution of equation \eqref{eq:g_n} at order $n$, satisfying the normalization \eqref{eq:ord-n-normalization}, are of the form $$g_n + \beta_n k^\perp, $$ where $g_n$ is uniquely determined and $\beta_n$ is arbitrary. 
        \end{itemize}
\end{lemma}

\begin{remark}
To prove Lemma \ref{lem:solution-induction}, \cite{Bus-Lla-23}, it is enough to prove that the average of the right hand side of equation \eqref{eq:g_n} is zero, that is, $\int_\torus R_n(\theta)d\theta = 0$. This is equivalent to prove that \begin{equation}
    \int_\torus k\cdot R_n(\theta) d\theta = 0\label{eq:k-Rn}
\end{equation} 
and \begin{equation}\label{eq:kperp-Rn}
    \int_\torus k^\perp \cdot R_n(\theta) d\theta = 0.
\end{equation}
A simple computation shows that \eqref{eq:k-Rn} is always satisfied, \cite{Bus-Lla-23}. Condition \eqref{eq:kperp-Rn} is attained by choosing $\beta_{n-1}$ appropriately. Note that \begin{equation}\label{eq:betan-conser}
    R_n(\theta) = \beta_{n-1} (n-1)! D^2V(\theta k + g_0) k^\perp  + S_n(\theta),
\end{equation}
where $S_n(\theta)$ depends on $g_0, \dots, g_{n-1}$. Thus, by the non-degeneracy condition \eqref{eq:non-degeneracy}, one can choose $\beta_{n-1}$ in such a way that \eqref{eq:kperp-Rn} is satisfied.
\end{remark}

\subsubsection{Dissipative case.}
In this case $\gamma = 1$ and the expansion $\mu_\eps = \sum \eps^n \mu_n $ also needs to be considered. To find a formal solution of \eqref{eq:inv-perpart} the following cohomology equations need to be solved:
\begin{align}
    \LL_\omega g_n(\theta)  &= R_n(\theta) + \mu_n, \qquad 0\leq n \leq 2,\label{eq:diss-ord-012}\\
    \LL_\omega g_3(\theta) &= R_3(\theta) - \omega k  + \mu_3, \label{eq:diss-ord-3} \\
    \LL_\omega g_n(\theta) &= R_n(\theta) + \mu_n - g_{n-3}(\theta) + g_{n-3}(\theta - \omega),\qquad n\geq 4, \label{eq:diss-ord-n}
\end{align}
where $\LL_\omega$ denotes the operator $\LL_\omega \varphi(\theta) := \varphi(\theta + \omega) - 2\varphi(\theta) + \varphi(\theta - \omega) $. These equations can be solved as follows:

For $n=0,1,2$ one can chose $\mu_n= 0$ and find  $g_n$ solving \eqref{eq:diss-ord-012}, and satisfying \eqref{eq:ord-n-normalization}, as in the conservative case. Note that, also in this case,  there exist at least two values of $\beta_0$ satisfying \eqref{eq:zero-ave-ord-1}. At order $n=3$ one chooses $\mu_3 = \omega k$ and proceeds also as in the conservative case.

For $n \geq 4$ one proceeds inductively. Assume that we have found $g_0, \dots, g_{n-1} + \beta_{n-1}k^\perp$ solutions of equations \eqref{eq:diss-ord-012}- \eqref{eq:diss-ord-n}, satisfying the normalization \eqref{eq:ord-n-normalization}, where $\beta_{n-1}$ is arbitrary. Equation \eqref{eq:diss-ord-n} of order $n$ can be solved as long as $\int_\torus \tilde R_n (\theta)d\theta = 0 $, where \begin{equation}\label{eq:zero-ave-diss}
    \tilde R_n (\theta) = R_n(\theta) +\mu_n - g_{n-3}(\theta) + g_{n-3}(\theta- \omega),
\end{equation} 
$R_n$ is computed as in \eqref{eq:R-n} and it depends on $g_0, g_1, \dots, g_{n-1}$. Note that $\int_\torus \tilde R_n(\theta)d\theta = 0$ is equivalent to \begin{equation}\label{eq:k-Rntilde}
    \int_\torus k\cdot \tilde R_n(\theta)d\theta = 0
\end{equation} 
and \begin{equation}\label{eq:kper-Rntilde}
    \int_\torus k^\perp\cdot \tilde R_n(\theta)d\theta = 0.
\end{equation}
Condition $\eqref{eq:k-Rntilde}$ is attained by adjusting the coefficiente $\mu_n$, that is, \begin{equation}
    \mu_n = -\left( \frac{1}{2\pi(k\cdot k)} \int_\torus k\cdot R_n(\theta) d\theta   \right)k. 
\end{equation}
To achieve condition \eqref{eq:kper-Rntilde} one needs to fix the constant $\beta_{n-1}$. Note that $R_n(g_0, g_1, \dots, g_{n-1} + \beta_{n-1}k^\perp )$ is the coefficient of order $\eps^{n-1}$ in the expansion of $\nabla V( \theta k + \sum_{i=0}^{n-1} g_i\eps^i + \beta_{n-1}\eps^{n-1}k^\perp )$, by computing the expansion one notes that  \begin{equation}\label{eq:betan-diss}
    \tilde R_n(\theta) = \beta_{n-1}(n-1)!D^2V(\theta k +g_0)k^\perp + S_n(\theta) +\mu_n -g_{n-3}(\theta) + g_{n-3}(\theta -\omega),
\end{equation}
where $S_n$ depends on $g_0, g_1, \dots, g_n, \beta_0, \dots, \beta_{n-2}$. Finally, if the non-degeneracy condition \eqref{eq:non-degeneracy} is fulfilled then $\beta_{n-1}$ can be fixed so that \eqref{eq:kper-Rntilde} is satisfied. Following the same procedure one can solve the equations to all orders and determine the formal solution $g_\eps, \mu_\eps$.

\begin{remark}
    We note that in both the conservative and the dissipative case the coefficients $R_n$ are computed using \eqref{eq:R-n}. However, due to the dependence on $g_0, g_1, \dots, g_{n-1}, \beta_0, \dots, \beta_{n-1}$, the values of $R_n$ will be different whether one solves $\eqref{eq:g_n}$ or \eqref{eq:diss-ord-012}-\eqref{eq:diss-ord-n}.
\end{remark}

\section{Numerical results.}\label{sec:num-results}
In this section we present the results of implementing the methods described in Section \ref{sec:lindstedt} and Appendix \ref{sec:appendix} to approximate numerically the domains of analyticity of the expansions of lower dimensional tori as functions of $\eps$.

First we have computed Lindstedt expansions \eqref{eq:g-mu-expan}, of the lower dimensional tori, as described in Section \ref{sec:lindstedt}. Once that the expansions were obtained, we fixed the value $\theta =1$ and computed Padé and Log-Padé approximants, see Appendix \ref{sec:appendix}. Finally, we computed the singularities of the Padé and Log-Padé approximants (the zeros of the denominators) to approximate the domains of analyticity of the expansions of the quasi-periodic orbits. 

\subsection{Comments on the implementations}
Equations \eqref{eq:g_n} and \eqref{eq:diss-ord-012}-\eqref{eq:diss-ord-n} have been solved using the Fourier representation of the functions, which is given in terms of their Fourier coefficients. The computation of the coefficients $R_n(\theta)$, in \eqref{eq:R-n}, is done by using a grid representation of the functions given by considering the values of the functions in an equally spaced grid. We switch between both representations by using the Discrete Fourier Transform.

It is proved in  \cite{Bus-Lla-23} that if $\nabla V$ is trigonometric polynomial of degree $J$, then $g_n$ is a trigonometric polynomial of degree at most $nJ$. This statement gives us a lower bound on the number of Fourier coefficients needed to obtain an accurate representation of the functions $g_n$ at order $n$.

We point out that we have used very high precisions to perform the computation due to the fact that the Lindstedt series are not expected to converge around $\eps=0$, see \cite{Jo-Lla-Zou-99, GallavottiG02}. Moreover, to prevent the accumulation of round-off errors we have used a filter that sets to zero the Fourier coefficients whose norm is smaller than $10^{-filter}$, where $filter\geq 100$.  We have used Pari/gp \cite{PARI2} to perform all the computations.

The explicit computation of the constants $\beta_n$ is done by taking the respective averages in \eqref{eq:betan-conser} and \eqref{eq:betan-diss} (in the conservative and the dissipative case), respectively, and equating to zero. All the computation were done using  $k=(1,0)$, $k^\perp = (0,1)$, $\beta_0 = 0$, and  $V(q) = -\cos(q_1) - \cos(q_2) - \cos(q_1 + q_2), q\in\torus^2$. It is not hard to verify that using this parameters \eqref{eq:zero-ave-ord-1} and \eqref{eq:non-degeneracy} are attained. 

Finally, we perform the computations using four different Diophantine frequencies $\omega= (\sqrt{5}-1)/2, (\sqrt{2})^{-1}, s^{-1}\approx 0.7548...$, and $\tau^{-1}\approx 0.5436$. Where $s$ is the real root of the polynomial $s^3 -s -1$  and $\tau$ is the real root of the polynomial $\tau^3 -\tau^2 -\tau -1$. All the four frequencies are well known Diophantine numbers. We note that with the notation introduced earlier $s^{-1}, \tau^{-1}\in \D(2)$ and $(\sqrt{5}-1)/2, (\sqrt{2})^{-1}\in \D(1)$, see \eqref{eq:diophantine}. So they have different Diophantine type.

\subsection{Domains of analyticity, conservative case.}
\begin{figure}[!ht]
    \centering
    \includegraphics[trim = 8cm 10cm 8cm 10cm, width=3.2truecm]{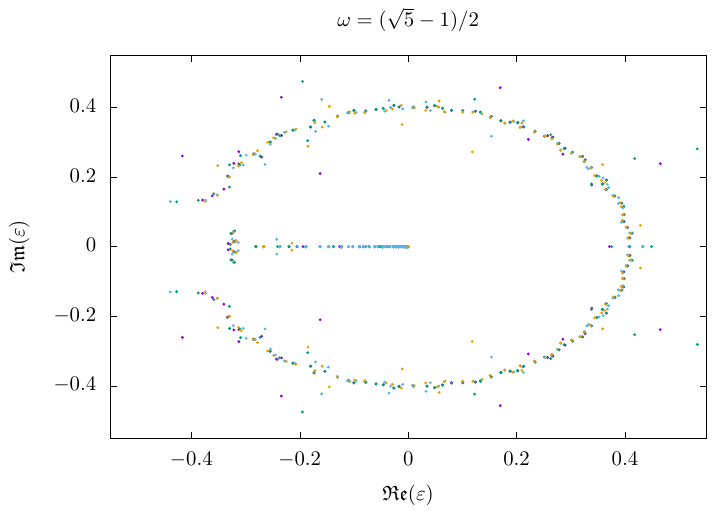}\hfil \hfil \hfil \hfil  
    \includegraphics[trim = 8cm 10cm 8cm 10cm, width=3.2truecm]{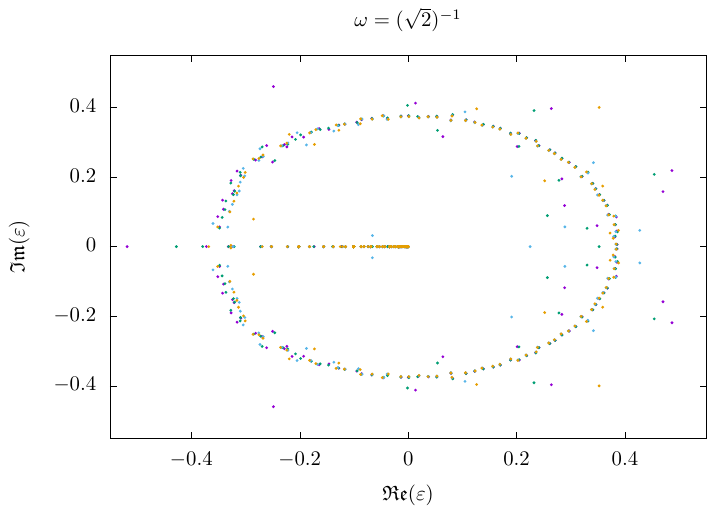}\\
    \bigskip 
    
    \includegraphics[trim = 8cm 10cm 8cm 10cm, width=3.2truecm]{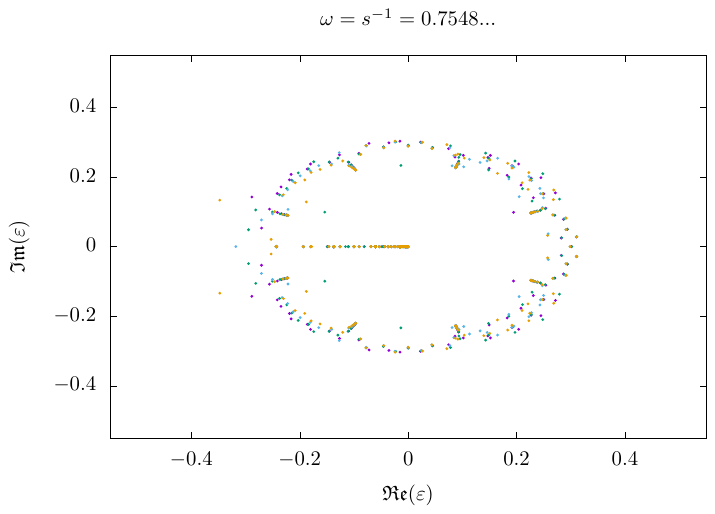}\hfil \hfil \hfil \hfil  
    \includegraphics[trim = 8cm 10cm 8cm 10cm, width=3.2truecm]{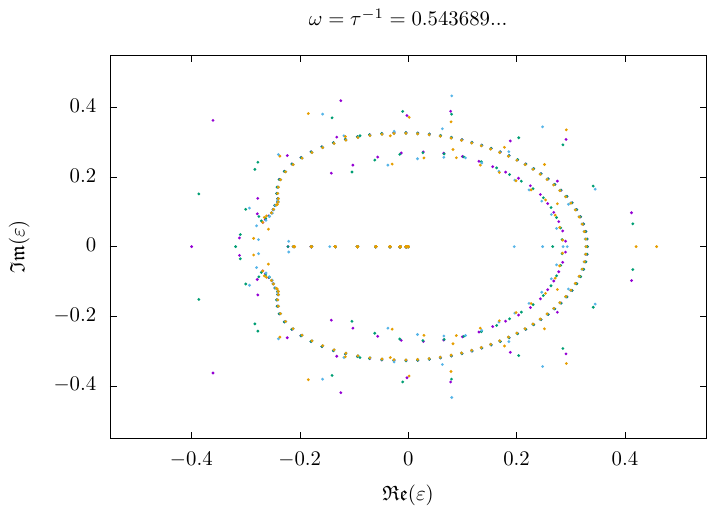}
    \caption{Approximation of the domains of analyticity of the Lindstedt series obtained by computing the poles of Padé and Log-Padé approximants for different Diophantine frequencies. Lindstedt series up to order $N=512$ were used. Conservative case.}
    \label{fig:domains-conservative}
\end{figure}

First we consider the case without dissipation, in Figure \ref{fig:domains-conservative} we plot the poles of the Padé and Log-Padé approximants of the Lindstedt series, the formal expansions were computed up to order $N=512$. We note that, for all the frequencies considered, there appear to be singularities in the negative real axis close to the origin, see Figure \ref{fig:domains-conservative}. As it is pointed out in \cite{Jo-Lla-Zou-99} when $\eps\in\real$, the case $\varepsilon <0$ corresponds to elliptic tori and,  moreover, resonances between the frequency of the torus and the normal frequencies happen when $\eps<0$. The case $\eps >0$ corresponds to whiskered tori. The approximation of the domains of analyticity depicted in the Figure \ref{fig:domains-conservative} support the conjectures in \cite{Jo-Lla-Zou-99, GallavottiG02} about the shape of the optimal domain of analyticity in $\eps$.

\subsection{Domains of analyticity, dissipative case.}
\begin{figure}[!ht]
    \centering
    \includegraphics[trim = 8cm 10cm 8cm 10cm, width=3.2truecm]{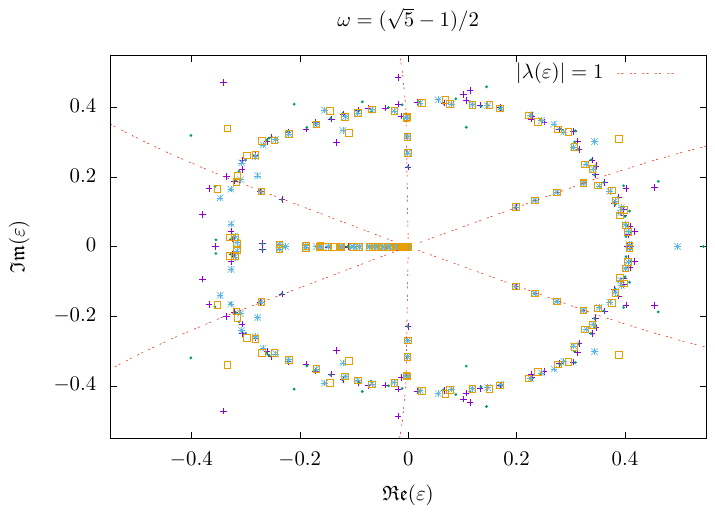}\hfil \hfil \hfil \hfil 
    \includegraphics[trim = 8cm 10cm 8cm 10cm, width=3.2truecm]{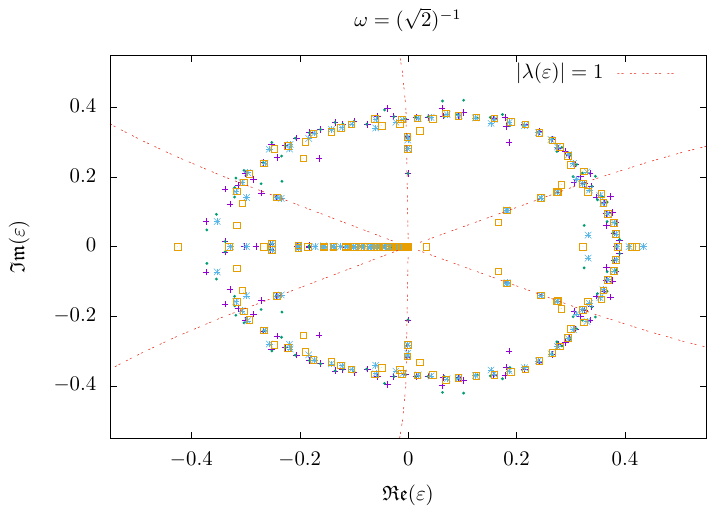}\\
    \bigskip 
    
    \includegraphics[trim = 8cm 10cm 8cm 10cm, width=3.2truecm]{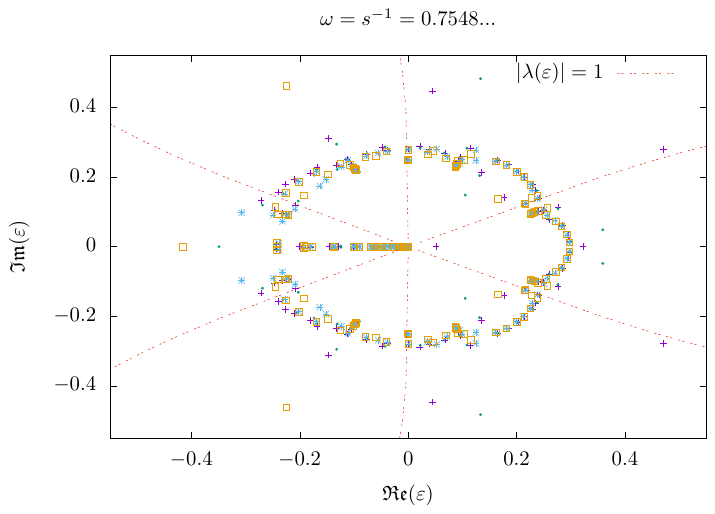}\hfil \hfil \hfil \hfil
    \includegraphics[trim = 8cm 10cm 8cm 10cm, width=3.2truecm]{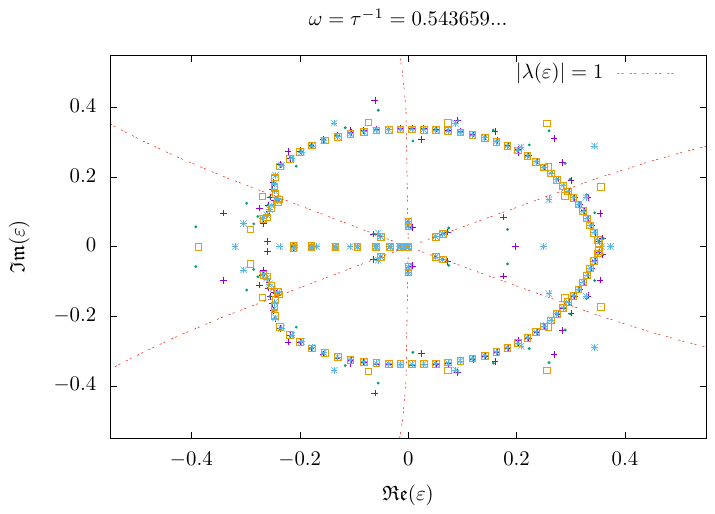}
    \caption{Approximation of the domains of analyticity of the Lindstedt series obtained by computing the poles of Padé and Log-Padé approximants for different Diophantine frequencies. Dissipative case, $\lambda(\varepsilon) = 1-\eps^3$.}
    \label{fig:domains-dissipative}
\end{figure}

The results for the dissipative case are summarized in Figure \ref{fig:domains-dissipative}, there we can observe that singularities in the negative real axis also appear, but in most of the cases there also appear some singularities over the lines described by the  equation $ |\lambda(\varepsilon)| = 1$ in the complex plane, with $\lambda(\varepsilon) = 1-\varepsilon^3 $. 

We recall that the proofs of existence of KAM tori and lower dimensional hyperbolic tori in \cite{CallejaCL17} and \cite{Cal-Cel-Lla-19-whisk}, respectively,  deal with small divisors coming from the conformal factor $\lambda(\eps)  = 1-\eps^3$. That is, small divisors appear for values of $\eps\in \complex$ such that $|\lambda(\eps)| = 1$.
We note that the use of  Lindstedt series usually involves less small divisors that the KAM methods. 

We remark that the main difference between \cite{Cal-Cel-Lla-19-whisk} and the present work is that in \cite{Cal-Cel-Lla-19-whisk} the authors assumed that for $\varepsilon=0$ the system has a whiskered tori of a given frequency and, as a consequence of this, the conjectured optimal domain of analyticity in \cite{Cal-Cel-Lla-19-whisk} only exhibits singularities for small values of $\eps\in \complex$ close to the curves described by the equations $|\lambda(\eps)|=1$. In the present work, when $\eps=0$, the map \eqref{eq:map} only has elliptic tori, the plots in Figure \ref{fig:domains-dissipative} seem to indicate that when $\eps<0$ some elliptic tori would be preserved despite of the presence of the dissipation. The plots also seem to suggest that when $\eps>0$ whiskered tori exist, as in the conservative case. 

\subsection{Gevrey character of Lindstedt series.}
The paper \cite{Bus-Lla-23} also contains results on the Gevrey character of the Lindstedt expansions of the lower dimensional tori. The approximations of the domains of analyticity found in the previous sections, see Figure \ref{fig:domains-conservative} and Figure \ref{fig:domains-dissipative} are in agreemetn with the Gevrey character of the series. We point out that the results in \cite{Bus-Lla-23} do not provide a sharp estimate of the Gevrey exponent so, as a byproduct of our computations, we give a numerical approximation of the Gevrey exponent of the series.

Given a fixed norm $\|\cdot \|$, we recall that a formal power series $\sum_{n=0}^\infty f_n \eps^n$ is called Gevrey if for any $n\in \nat$ there exist $A, R, \sigma >0$ such that \begin{equation}\label{eq:gevrey-def}
    \|f_n\| \leq AR^n n^{\sigma n}.
\end{equation}
The constant $\sigma$ in \eqref{eq:gevrey-def} is called the Gevrey exponent. Note that in this work we consider formal power series of the form $\sum g_n\eps^n$ where $g_n:\real\longrightarrow \real^2$ are periodic functions. Given a periodic function, $g$, with Fourier expansion $g(\theta) = \sum_{\ell\in \integer} \hat g_\ell e^{i\theta \ell}$  the norm used in \cite{Bus-Lla-23} to prove the Gevrey character of the series was the following \begin{equation}
    \|g\|^2_{\rho, r} := \sum_{\ell\in \integer} |\hat g_\ell|^2 e^{2|\ell|\rho} \left( 1 + |\ell|\right)^r
\end{equation} 
where $\rho\geq 0$ and $r\in \integer_+$. We recall that all the functions in our expansions are trigonometric polynomials, see \cite{Bus-Lla-23}, so we will always be computing finite sums.

\begin{figure}[!ht]
    \centering
    \includegraphics[trim = 8cm 10cm 8cm 10cm, width=3.2truecm]{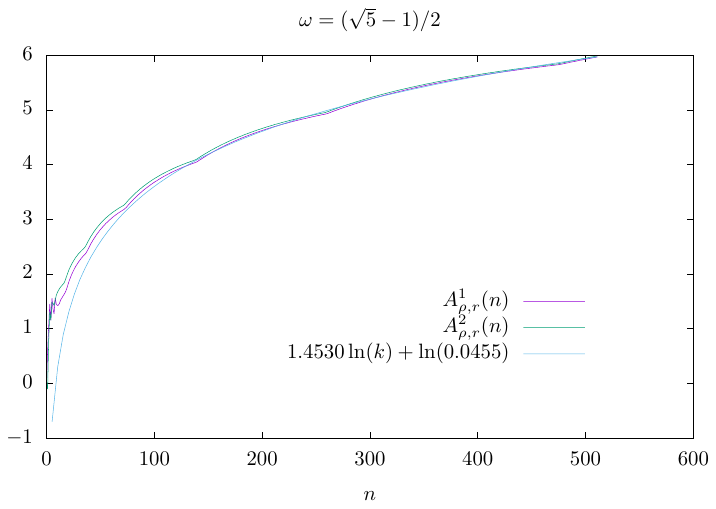}\hfil \hfil \hfil \hfil 
    \includegraphics[trim = 8cm 10cm 8cm 10cm, width=3.2truecm]{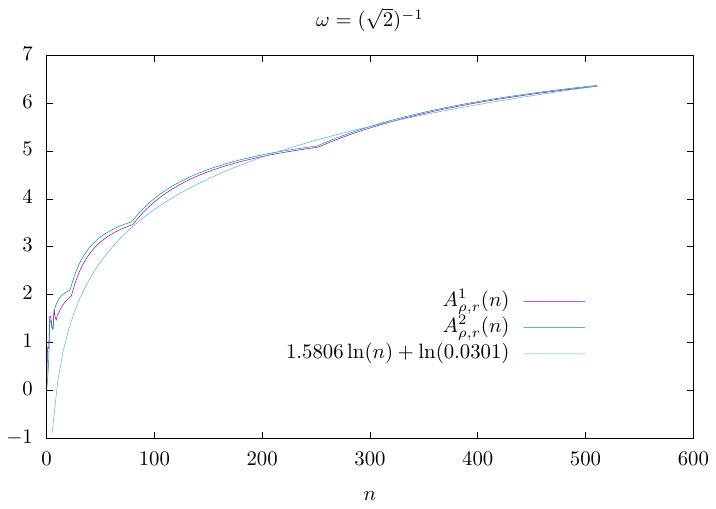}\\
    \bigskip 
    
    \includegraphics[trim = 8cm 10cm 8cm 10cm, width=3.2truecm]{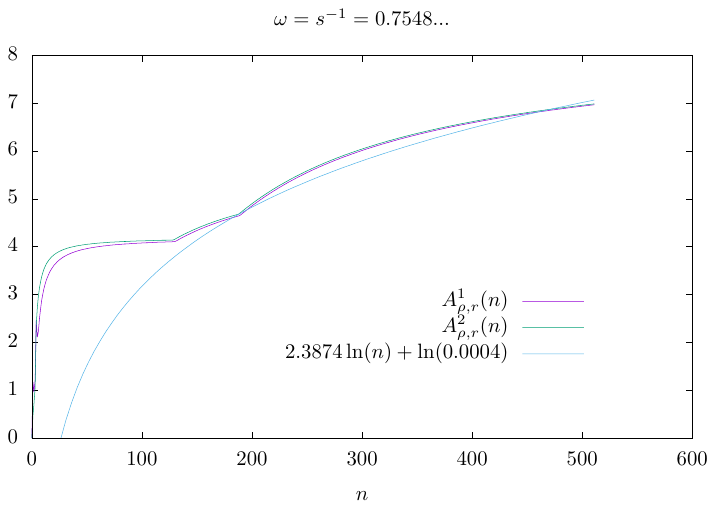}\hfil \hfil \hfil \hfil
    \includegraphics[trim = 8cm 10cm 8cm 10cm, width=3.2truecm]{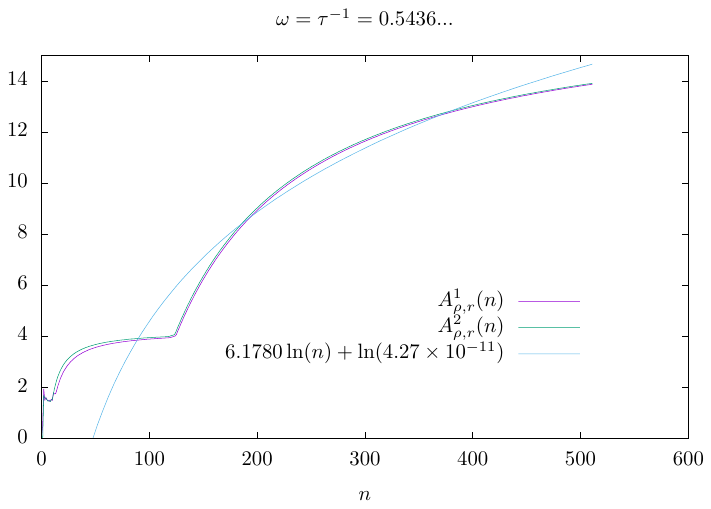}
    \caption{Plots of $A_{\rho,r}^i(n) =\frac{1}{n} \ln\left(\|g_n^i\|_{\rho,r} \right) $  and it respective numerical fits, where $g_n = (g_n^1, g_n^2)$.  Different values of the frequency $\omega$ were used and $\rho = 0.1$, and $r=4$. Conservative case. }
    \label{fig:norms-fits-conservative}
\end{figure}

We follow the approach of \cite{BustamanteC19, BustamanteC21} to estimate numerically the Gevrey character of the Lindstedt series we have computed. That is, given the power series $\sum g_n \eps^n$ we define \begin{equation}\label{eq:a-rho-r}
    A_{\rho,r} (n) :=\frac{1}{n} \ln\left(\|g_n\|_{\rho,r} \right) . 
\end{equation}
Note that if the coefficients of $\sum g_n\eps^n$ satisfy \eqref{eq:gevrey-def} then \begin{equation}
A_{\rho,r} (n) \sim  \sigma \ln(n) + \ln(R).     
\end{equation} 
Thus, a logarithmic growth of $A_{\rho,r}$ would give numerical evidence of the Gevrey character of the series. To estimate numerically the Gevrey exponent, $\sigma$, we have computed  $A_{\rho, r}(n)$ and then fit a logarithmic function of the form $\sigma \ln(n) + \ln(R)$.   The results for the conservative case are summarized in Table \ref{tab:fits-conservative} and Figure \ref{fig:norms-fits-conservative}. We only show results for $\rho =0.1$ and $r=4$, results for different values of these constants are very similar.
The results in Table \ref{tab:fits-conservative} seem to indicate that the Gevrey exponent depends on the frequency. 

\begin{table}[ht]
\begin{tabular}{ |p{2cm}||p{2cm}|p{2cm}|}
 \hline
 \multicolumn{3}{|c|}{Conservative case} \\
 \hline
 & $R$ & $\sigma$ \\
 \hline
 $ \omega = \frac{\sqrt{5}-1}{2} $ &  0.045441 & 1.453009\\ 
 $ \omega = \sqrt{2}  $ &   0.030117  & 1.580654\\ 
 $\omega = s^{-1} $ &  0.000472 &  2.387420 \\
 $\omega = \tau^{-1} $ & 4.27e-11  & 6.178068 \\
 \hline
\end{tabular}

\caption{Numerical fit of a function $\log(R) +\sigma\log(n)$ to the data  $A_{\rho,r}(k)$ for $\rho =0.1$, $r=4$ and  different values of the frequency $\omega$. Computations were done using $2^{11}$ Fourier coefficients and at least 1500 digits of precision. The numerical fit was made in for $150\leq n\leq 450$.}
\label{tab:fits-conservative}
\end{table}

\begin{figure}[!ht]
    \centering
    \includegraphics[trim = 8cm 10cm 8cm 10cm, width=3.4truecm]{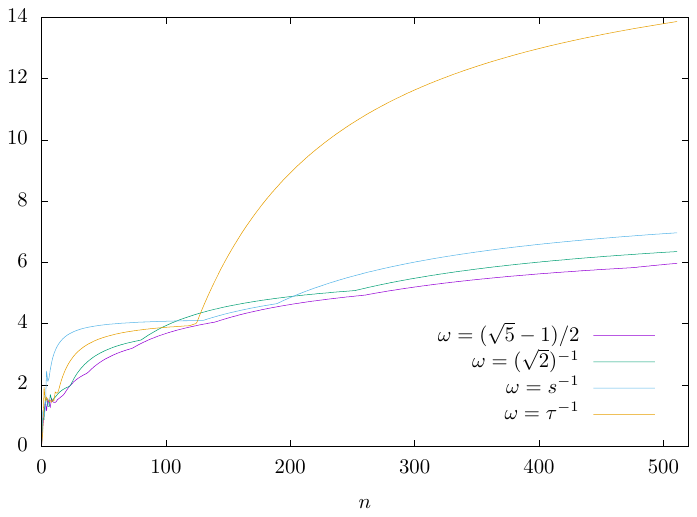}
    \caption{Plots of $A_{\rho, r}^1(n) = \frac{1}{n} \ln\left(\|g_n^1\|_{\rho,r} \right)$ for different values of $\omega$, $\rho=0.1$, $r=4$. Dissipative case.}
    \label{fig:norms-comparison-dissipative}
\end{figure}

The growth of the coefficients for the dissipative case is very similar as in the conservative case. At first sight the coefficients $A_{\rho,r}(n)$ look the same as in the conservative case, see Figure \ref{fig:norms-comparison-dissipative}. However, the quantities $A_{\rho, r}(n)$ are not the same  in the conservative and dissipative case but the digits up to order $10^{-6}$ coincide,  that is why the graphs and the fits look the same for the conservative and dissipative case, see Table \ref{tab:comparison-coefficients}. 

We recall that in \cite{BustamanteL22} was showed that, in the dissipative case, the formal expansions of KAM tori of a given frequency, $\omega$, are Gevrey and, moreover, an upper bound of the Gevrey exponent, $\sigma$, was also given, that is, $\sigma\leq \ell/3$, where $\omega\in \D(\ell)$ as in \eqref{eq:diophantine}. Table \ref{tab:fits-conservative}, which is the same for the conservative and dissipative case, give numerical evidence that the Gevrey exponent of expansions of lower dimensional tori is larger than the exponent for expansions of KAM tori. Note that $(\sqrt{5}-1)/2, (\sqrt{2})^{-1} \in \D(1)$ and $s^{-1},\tau^{-1}\in \D(2)$.

\begin{table}[ht]

\begin{tabular}{ |p{2cm}||p{2.2cm}|p{2.2cm}|}
 \hline
 \multicolumn{3}{|c|}{$A_{\rho,r}(n)$} \\
 \hline
 & Conservative  & Dissipative \\
 \hline
 $ n=20 $ &  2.8503536756  & 2.8503522703\\ 
 $ n=40 $ &  3.5283833570  & 3.5283818284 \\ 
 $n=60$ &  3.7543932509 &  3.7543916811 \\
 $n=80 $ & 3.8673981978  & 3.8673966074 \\
 $n=100$ & 3.9352011659 & 3.9351995632 \\
 $n=120$  &  4.0064464018  &  4.0064463738 \\
 \hline
\end{tabular}

\caption{Comparison of the coefficients $A_{\rho,r}(n)$ for the conservative and dissipative cases when $\omega = \tau^{-1}$, $\rho=0.1$, and $r=2$. }
\label{tab:comparison-coefficients}
\end{table}

\subsection{Validation of the results.}
\begin{figure}[!ht]
    \centering
    \includegraphics[trim = 8cm 10cm 8cm 10cm, width=3.2truecm]{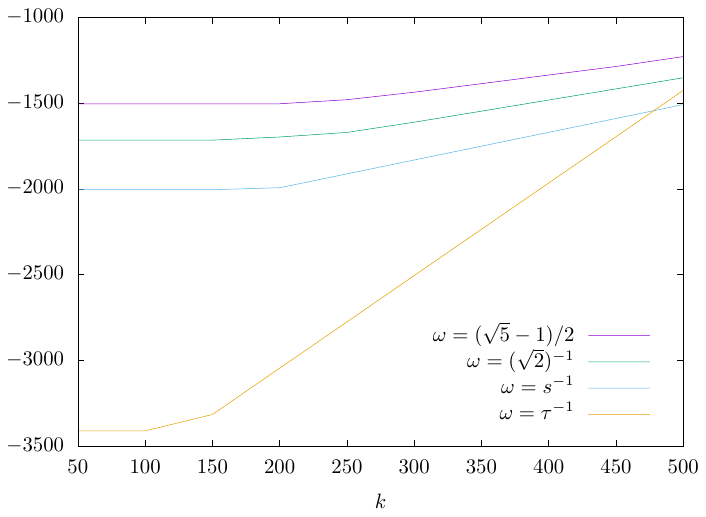}\hfil \hfil \hfil \hfil 
    \includegraphics[trim = 8cm 10cm 8cm 10cm, width=3.2truecm]{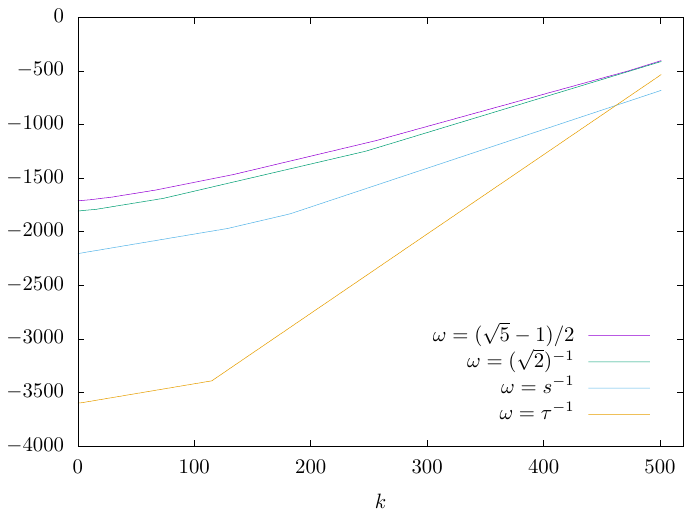}\\
    
    \caption{Left panel: Conservative case. Error of the invariance equation\eqref{eq:inv-perpart} up to order $n$ in logarithmic scale, that is  $\ln\left(\| \LL_\omega \sum_{i=0}^n g_i\eps^i - \eps \nabla V(\theta k + \sum_{i=0}^n g_i(\theta)\eps^i) \|_\infty\right)$, with $\eps=10^{-2}$.  Right panel: Dissipative case. Error of the cohomology equation \eqref{eq:diss-ord-n} at order $n$ in logarithmic scale, that is $\ln\left(\|\LL g_n - R_n-\mu_n +g_{n-3}-g_{n-3}\circ T_{-\omega}\|_\infty  \right)$. }
    \label{fig:errors-cons-diss} 
\end{figure}

All the computations were done with at least 1500 digits of precision. High precision is needed due to the rapid growth of the coefficients of the Lindstedt series $g_n$. Figure \ref{fig:errors-cons-diss}, left panel, shows plots of the norm of the error of the invariance equation \eqref{eq:inv-perpart} in the conservative case when we evaluate \eqref{eq:inv-perpart} in the partial sums $\sum_{i=0}^n g_i(\theta)\eps^i$, with $\eps =10^{-2}$. The right panel of Figure \ref{fig:errors-cons-diss}, shows plots of the norm of the cohomology equations \eqref{eq:diss-ord-n} at order $n$ for the dissipative case. All the computation were done using Pari/GP, \cite{PARI2}.

\section{Conclusion}
We considered a Froeschlé map in its conservative and weakly dissipative versions. We have implemented numerically an algorithm to compute Lindstedts series of lower dimensional tori and use them to approximate their domains of analyticity as functions of a complex parameter of perturbation $\eps\in \complex$. 

Our results support conjectures in the literature \cite{Jo-Lla-Zou-99, GallavottiG02} about the shape of the optimal domain of analyticity of lower dimensional tori in the conservative case. When the perturbation also introduces a small dissipation of the form $\lambda(\eps) = 1-\eps^3$, the domains of analyticity seem to exhibit more singularities, which have appeared in the literature before \cite{CallejaCL17, Cal-Cel-Lla-19-whisk, BustamanteC19}. Also in the dissipative case, our results suggest the persistance of some elliptic tori when $\eps<0$ and
the existence of whiskered tori when $\eps>0$.

We also studied numerically the size of the Gevrey exponent of the Lindstedt series for lower dimensional tori. Our results indicate that in both the conservative and dissipative case the Gevrey exponent depends on the Diophantine type of the frequency. We also found numerical evidence that the Gevrey exponent of the expansions of lower dimensional tori is larger than the exponent for expansions of KAM tori in the weakly dissipative case. However, our computations do not allow us to distinguish any difference between the exponents in the dissipative and conservative case.

\begin{appendix}
\section{Padé and Log-Padé approximants}\label{sec:appendix}
We briefly describe the construction of a Padé approximant. This method is well know and can be found in several places in the literature, see for example \cite{baker1996pade}. 

A Padé approximant  of order $[m/n]$ of a power series $f_\eps = \sum_{i=0}^\infty f_i\eps^i$ is a rational function $M(\eps)/N(\eps)$, with $N(0)=1$, which agrees with $f_\eps$ to the highest possible order. That is \begin{equation}\label{eq:pade1}
    f_\eps - \frac{M(\eps)}{N(\eps)} = \mathcal{O}\left(|\eps|^{m+n+1}\right),
\end{equation}
where $M(\eps)$ and $N(\eps)$ are polynomials of degrees $m$ and $n$, respectively.

Denote $M(\eps) = \sum_{i=0}^m M_i\eps^i$ and $N(\eps) = \sum_{i=0}^n N_i\eps^i$, then the coefficients $M_i,N_i$ can be found solving the following linear system of equations\begin{align}
    g_i + \sum_{i=1}^j g_{j-i}N_i &= M_j\quad 0\leq j\leq m \\
    g_i + \sum_{i=1}^n g_{j-i}N_i &= 0 \quad m < j \leq m+n 
\end{align} 

Padé approximants are better suited to approximate functions with poles rather than functions with branch points. The works  \cite{Lla-Tom-94, Lla-Tom-95} indicate that the boundaries of the domains of analyticity of quasi-periodic orbits could be described as an accumulation of branch points. Then, a logarithmic Padé approximation would be better suited to approximate singularities that are branch points. 

The computation of a logarithmic Padé approximation could be done as follows, see \cite{Lla-Tom-95}. If a function $f(\eps)$ has a branch point singularity of order $\ell\in \integer$ at $\eps = \eps_0$, then $f(\eps) = A(\eps - \eps_0)^{-1/\ell} + g(\eps) $ with $g$ and analytic function at $\eps_0$. For $\eps $ close to $\eps_0$ one has \begin{equation}\label{eq:log-pad-func}
    F(\eps):= \frac{d}{d\eps}\ln(f(\eps)) = \frac{f'(\eps)}{f(\eps)}\approx -\frac{\frac{1}{\ell}}{\eps-\eps_0},
\end{equation}   
so one expects that a Padé approximant of $F$ to exhibit a pole at $\eps=\eps_0$ with residue $-1/\ell$. We call a \textit{Log-pade approximation of $f$} to a Padé approximant of the function $F(\eps)$ defined in \eqref{eq:log-pad-func}.

\end{appendix}

\bibliographystyle{alpha}
\bibliography{bibliography}

\end{document}